\documentclass[12pt]{article}
\usepackage{amsfonts,amsmath,amsxtra}

\newtheorem{theorem}{Theorem}

\newtheorem{lemma}{Lemma}

\newenvironment{definition}
{\smallskip\noindent{\bf Definition\/}:}{\smallskip\par}
\newenvironment{proposition}
{\smallskip\noindent{\bf Proposition\/}.}{\smallskip\par}

\newenvironment{proof}
{\noindent{\bf Proof\/}.}{{ $\Box$}\smallskip\par}

\title{Motivic integrals and functional equations.}
\author{E. Gorsky\footnote{Supported by the grant NSh-4719.2006.1. \newline AMS 2000 Subj. Class.:32S45, 28B10. }}
\date{}

\begin{document}

\maketitle

\begin{abstract}

A functional equation for the motivic integral corresponding to the Milnor number of an arc is derived using the Denef-Loeser formula for the
change of variables. Its solution is a function of five auxiliary parameters, it is unique up to multiplication by a constant and there is a simple recursive algorithm to
find its coefficients.
The method is universal enough and gives,
for example, equations  for the integral corresponding to the intersection number
over the space of pairs of arcs and over the space of unordered tuples of arcs.
\end{abstract}

\section{ Introduction}

Motivic integration, introduced by M. Kontsevich, is a powerful tool for exploring the space of formal arcs on
a given variety.  Motivic integrals provide the generating series
for motivic measures of level sets of some arc invariants. There are examples  of such integrals 
which can be calculated explicitly
(see e.g. \cite{dl}). In more general
situation the values of such integrals are unknown, however some of them  satisfy  functional equations if some
auxiliary variables are introduced.

In this paper a functional equation for the motivic integral which gives the generating series corresponding to the Milnor number of a plane curve is derived using the Denef-Loeser formula for the
change of variables (\cite{dl}). Its solution is unique up to  multiplication by a constant and there is a simple algorithm to
express its coefficients via the initial ones.
For example it implies partial differential equations
for the solution. This equation gives a method
to compute the motivic measure (and, consequently,
the Hodge-Deligne polynomial) of the stratum
 $\{\mu=const\}$ in the space of the plane curves.
Some examples are considered in Section 4.  

A similar idea gives some other equations,
for example, an equation for the integral corresponding to the intersection number over the 
space of pairs of arcs. Moreover, using the notion of the power structure
on the Grothendieck ring (\cite{powers}) we introduce a motivic measure on the space of unordered tuples of arcs. A curious equation for the integral corresponding to the intersection number in this case is derived as well.

Some generating series with coefficients in the Grothendieck ring of varieties
 $K_0(Var_{\mathbb{C}})$
(or in the Grothendieck ring of Chow motives) satisfy functional
equations similar to the functional equation for Hasse-Weil zeta function. These equations, obtained by M. Kapranov (\cite{kapr}) and
F. Heinloth (\cite{hein}), follow from the duality theory on curves
and Abelian varieties. Notice that  they have 
 origin different from our approach.

\section{ Motivic measure}

Let $\mathcal{L}=\mathcal{L}_{\mathbb{C}^2,0}$ be the space
of arcs at the origin on the plane. 
It is the set of pairs $(x(t),y(t))$ of formal power series 
(without degree 0 terms).

Let $K_0(Var_{\mathbb{C}})$ be the Grothendieck ring of quasiprojective
complex algebraic varieties. It is generated by the isomorphism classes
of complex quasiprojective algebraic varieties modulo relations
 $[X]=[Y]+[X\setminus Y],$ where $Y$ is a Zariski closed subset of  $X$. 
Multiplication is given by the formula $[X]\cdot [Y]=[X\times Y].$
Let $\mathbb{L}\in K_0(Var_{\mathbb{C}})$ be the class of the complex line.

Consider the ring $K_0(Var_{\mathbb{C}})[\mathbb{L}^{-1}]$
with the following filtration: $F_k$ is generated by the elements
of type $[X]\cdot [\mathbb{L}^{-n}]$ with $n-\dim X\ge k$. Let $\mathcal{M}$
be the completion of $K_0(Var_{\mathbb{C}})[\mathbb{L}^{-1}]$ corresponding
to this filtration.

On an algebra of subsets of the space $\mathcal{L}$ J. Denef and F. Loeser (\cite{dl}) 
(after M. Kontsevich) have constructed a measure $\chi_g$ with values in the ring  $\mathcal{M}$. According to this measure, one can naturally define the (motivic) integral for simple functions on $\mathcal{L}$(\cite{dl}). 

We will use the simple functions
$v_{x}=Ord_{0}x(t), v_y=Ord_{0}y(t)$ and $v=\min\{v_x,v_y\}$,
defined for an arc $\gamma(t)=(x(t),y(t)).$

Let $h:Y\rightarrow X$ be a proper birational morphism
of smooth manifolds of  dimension $d$
and $J=h^{*}K_{X}-K_{Y}$ be the relative canonical divisor
on  $Y$ (locally it is defined by the Jacobi determinant).
It defines a function $ord_{J}$ on the space of arcs on $Y$ --
the intersection number between the arc and the divisor.
Then one has the following change of variables formula in the motivic 
integral:
  
\begin{theorem}(\cite{dl})
Let $A$ be a measurable subset in the space of arcs on
$X$, let $\alpha$ be a simple function. Then
$$\int_{A}\alpha d\chi_g=\int_{h^{-1}(A)}(h^{*}\alpha)\mathbb{L}^{-ord_{J}}d\chi_g.$$
\end{theorem}

If $h$ is a blow-up of the origin in the plane, the relative
canonical divisor coincides with the exceptional line, so the
function $ord_{J}$ coincides with the intersection number with this line.

\section{ Functional equation for Milnor number}

The Milnor number of the plane curve given by the equation $\{f=0\}$
can be defined as the codimension of the ideal generated
by the partial derivatives ${\partial f\over\partial x},{\partial f\over\partial y}$.

We shall use the following statement.

\begin{lemma}
Suppose that upon after blowing up the origin the Milnor number
of an irreducible curve is equal to
$\mu$, and the intersection number with the exceptional divisor
is equal to  $p$. Then the Milnor number of the initial curve is equal to
$$\mu+p(p-1).$$
\end{lemma}


Let 
$$I(t,a,b,c,d,f)=\int_{\mathcal{L}}t^{\mu}a^{v_x}b^{v_y}c^{v_x^2}d^{v_xv_y}f^{v_y^2}d\chi_{g}.$$

\begin{theorem}
This function satisfies the functional equation
$$I(t,a,b,c,d,f)=I(t, t^{-1}ab\mathbb{L}^{-1}, b, tcdf, df^2, f)+
I(t, t^{-1}ab\mathbb{L}^{-1}, a, tcdf, dc^2, c)+$$
$$+I(t, t^{-1}ab\mathbb{L}^{-1}, 1, tcdf, 1, 1)\cdot(\mathbb{L}-1).$$
\end{theorem}

\begin{proof}
Let
$$A(t,a,b,c,d,f)=\int_{\{v_y>v_x\}}t^{\mu}a^{v_x}b^{v_y}c^{v_x^2}d^{v_xv_y}f^{v_y^2}d\chi_{g},$$
and note that
$$\int_{\{v_x>v_y\}}t^{\mu}a^{v_x}b^{v_y}c^{v_x^2}d^{v_xv_y}f^{v_y^2}d\chi_{g}=A(t,b,a,f,d,c).$$
Let us compute the analogous integral over $\{v_x=v_y\}$.
For $v_x=v_y$ $$y=\lambda x+\widetilde{y},\lambda\neq 0, v(\widetilde{y})>v(y).$$
For $\lambda$ fixed $\mu(x(t),y(t))=\mu(x(t),\widetilde{y}(t)),$
hence
$$\int_{\{v_x=v_y\}}t^{\mu}a^{v_x}b^{v_y}c^{v_x^2}d^{v_xv_y}f^{v_y^2}d\chi_{g}=
\int_{\{v_x=v_y\}}t^{\mu}(ab)^{v_x}(cdf)^{v_x^2}d\chi_{g}=$$
$$(\mathbb{L}-1)\int_{\{v_x<v_{\widetilde{y}}\}}t^{\mu}(ab)^{v_x}(cdf)^{v_x^2}d\chi_{g}=
A(t,ab,1,cdf,1,1)\cdot(\mathbb{L}-1),$$
and consequently,
$$I(t,a,b,c,d,f)=A(t,a,b,c,d,f)+A(t,b,a,f,d,c)+A(t,ab,1,cdf,1,1,1)\cdot(\mathbb{L}-1).$$

Let us blow-up the origin. If $v_y>v_x$, then $y(t)=x(t)\theta(t),\theta(0)=0$,
and therefore the corresponding modifications of the curves pass through 
the fixed point $p_0$ of the exceptional divisor related to the $x$-axis. Thus
$$\sigma^{*}\mu=\mu+v_x(v_x-1), \sigma^{*}v_x=v_x, \sigma^{*}v_x^2=v_x^2,$$
$$\sigma^{*}v_y=v_y+v_{\theta}, \sigma^{*}v_xv_y=v_x^2+v_xv_{\theta}, \sigma^{*}v_y^2=v_x^2+2v_xv_{\theta}+v_{\theta}^2.$$
Using the Denef-Loeser change of variables formula, we obtain:
$$A(t,a,b,c,d,f)=\int_{\{v_y>v_x\}}t^{\mu}a^{v_x}b^{v_y}c^{v_x^2}d^{v_xv_y}f^{v_y^2}d\chi_{g}=$$
$$\int_{\mathcal{L}}t^{\mu+v_x(v_x-1)}a^{v_x}b^{v_x+v_{\theta}}c^{v_x^2}d^{v_x^2+v_xv_{\theta}}f^{v_x^2+2v_xv_{\theta}+v_{\theta}^2}\mathbb{L}^{-v_x}d\chi_{g}=$$
$$\int_{\mathcal{L}}t^{\mu}(t^{-1}ab\mathbb{L}^{-1})^{v_x}b^{v_{\theta}}(tcdf)^{v_x^2}(df^2)^{v_xv_{\theta}}f^{v_{\theta}^2}d\chi_{g}=
I(t, t^{-1}ab\mathbb{L}^{-1}, b, tcdf, df^2, f),$$
therefore 

\begin{equation}
\label{mu}
I(t,a,b,c,d,f)=I(t, t^{-1}ab\mathbb{L}^{-1}, b, tcdf, df^2, f)+
\end{equation}
$$I(t, t^{-1}ab\mathbb{L}^{-1}, a, tcdf, dc^2, c)+
+I(t, t^{-1}ab\mathbb{L}^{-1}, 1, tcdf, 1, 1)\cdot(\mathbb{L}-1).$$
\end{proof}

It is clear that

\begin{equation}
\label{sym}
I(t,a,b,c,d,f)=I(t,b,a,f,d,c)
\end{equation}

and differentiating under the integral we get

\begin{equation}
\label{pde}
-c{\partial I\over\partial c}=[-a{\partial \over\partial a}]^2I,-d{\partial I\over\partial d}=[a{\partial \over\partial a}]\circ b{\partial \over\partial b}I,\\
-f{\partial I\over\partial f}=[-b{\partial \over\partial b}]^2I.
\end{equation}

\begin{theorem}
A function $I(t,a , b, c, d, f)$, divisible by $abcdf$, 
 satisfying the functional equation
(\ref{mu}) and the symmetry condition (\ref{sym}),  is unique
up to multiplication by a constant.
\end{theorem}

Before proving Theorem 3 let us consider  the following simpler example of an analogous functional
equation for the motivic integral. Let
$$f(a,b)=\int_{\mathcal{L}}a^{v_x}b^{v_y}d\chi_{g}={ab\mathbb{L}^{-2}(\mathbb{L}-1)^2\over (1-a\mathbb{L}^{-1})(1-b\mathbb{L}^{-1})}.$$
Similarly to the calculations above using the change of variables formula
one can obtain the functional equation
\begin{equation}
\label{3,5}
f(a,b)=f(ab\mathbb{L}^{-1},a)+f(ab\mathbb{L}^{-1},b)+f(ab\mathbb{L}^{-1},1)\cdot(\mathbb{L}-1).
\end{equation}
Let us describe its solutions. Note that from the definition of $f(a,b)$
it follows  that 
$f(a,b)=f(b,a), f(0,b)=0.$ 

Let $f(a,b)=\sum_{i,j}f_{ij}a^ib^j.$ Then the functional equation
(\ref{3,5}) can be rewritten in the form
$$\sum_{i,j}f_{ij}a^ib^j=\sum_{i,j}f_{ij}\mathbb{L}^{-i}a^{i+j}b^{i}+
\sum_{i,j}f_{ij}\mathbb{L}^{-i}a^{i}b^{i+j}+(\mathbb{L}-1)\sum_{i,j}f_{ij}\mathbb{L}^{-i}a^{i}b^{i}=$$
$$\sum_{i\ge j}f_{j,i-j}\mathbb{L}^{-j}a^{i}b^{j}+
\sum_{i\le j}f_{i,j-i}\mathbb{L}^{-i}a^{i}b^{j}+(\mathbb{L}-1)\sum_{i,j}f_{ij}\mathbb{L}^{-i}a^{i}b^{i}.$$
Using the relations $f_{ij}=f_{ji}, f_{i0}=0,$ we  obtain 
the system of recurrence relations on the coefficients:
\begin{equation}
\begin{cases}
f_{ij}=\mathbb{L}^{-j}f_{i-j,j}, i>j\\
f_{ii}=\mathbb{L}^{-i}(\mathbb{L}-1)\sum_{j=1}^{\infty}f_{ij}.\\
\end{cases}
\end{equation}

Below we will use a more general system of equations:

\begin{equation}
\label{recur}
\begin{cases}
f_{ij}=\varepsilon_{j}f_{i-j,j}, i>j,\\
f_{ii}=C\varepsilon_{i}\sum_{j=1}^{\infty}f_{ij}.\\
\end{cases}
\end{equation}

\begin{lemma}
Let $1-\varepsilon_i-C\varepsilon_i\neq 0,C\neq 0, \varepsilon_{i}\neq 0,\varepsilon_{i}\neq 1$ for all $i$.
Then the system (\ref{recur}) has a non-zero solution, which is
defined uniquely up to multiplication by a constant, such that
 $f_{i,j}=f_{j,i}$ and $f_{i,0}=0.$
\end{lemma}

\begin{proof}
From the first equation of (\ref{recur}) one has
$$\sum_{j>i}f_{ij}=\varepsilon_{i}\sum_{j>i}f_{i,j-i}=\varepsilon_{i}\sum_{j>0}f_{i,j}.$$
Moreover, $f_{ii}=C\varepsilon_{i}\sum_{j>0}f_{ij},$ so
$f_{ii}=C\sum_{j>i}f_{ij}, $ thus 
$$\sum_{0<j<i}f_{ij}=\sum_{j>0}f_{ij}-f_{ii}-\sum_{j>i}f_{ij}=
f_{ii}({1\over C\varepsilon_{i}}-1-{1\over C}),$$
so $$f_{ii}={C\varepsilon_{i}\over 1-\varepsilon_{i}-C\varepsilon_{i}}\sum_{0<j<i}f_{ij}.$$

Let $f_{11}$ be an arbitrary non-zero number. Let us compute
$f_{ij}$. If $i\neq j$, we can use the first equation of
(\ref{recur}), and if $i=j$  we can use the previous equality.
In any case $f_{ij}$ will be expressed via
$f_{k,l}$ with $k+l<i+j$, so this process will stop
and $f_{ij}$ will be expressed via
$f_{11}$. Therefore the solution is unique.

It is easy to prove that the described algorithm defines the solution to (\ref{recur}).

\end{proof}

Let us prove Theorem 3.
Consider the equation (\ref{mu}). Let 
$$I(t,a,b,c,d,f)=\sum_{k_1,k_2,k_3,k_4,k_5}g_{k_1,k_2,k_3,k_4,k_5}a^{k_1}b^{k_2}c^{k_3}d^{k_4}f^{k_5},$$
then 

\centerline{$\sum_{\underline{k}}g_{\underline{k}}a^{k_1}b^{k_2}c^{k_3}d^{k_4}f^{k_5}=\sum_{k_1,k_2,k_3,k_4,k_5}g_{k_1,k_2,k_3,k_4,k_5}(t^{-1}\mathbb{L}^{-1}ab)^{k_1}b^{k_2}(tcdf)^{k_3}(df^2)^{k_4}f^{k_5}+$}
$$\sum_{k_1,k_2,k_3,k_4,k_5}g_{k_1,k_2,k_3,k_4,k_5}(t^{-1}\mathbb{L}^{-1}ab)^{k_1}a^{k_2}(tcdf)^{k_3}(dc^2)^{k_4}c^{k_5}+$$
$$(\mathbb{L}-1)\cdot\sum_{k_1,k_2,k_3,k_4,k_5}g_{k_1,k_2,k_3,k_4,k_5}(t^{-1}\mathbb{L}^{-1}ab)^{k_1}(tcdf)^{k_3}.$$

We obtain the following system of equations
\begin{equation}
\label{sys}
\begin{cases}
g_{k_1,k_2,k_3,k_4,k_5}=t^{k_3-k_1}\mathbb{L}^{-k_1}g_{k_1,k_2-k_1,k_3,k_4-k_3,k_5-2k_4+k_3},\mbox{\rm if}\quad k_2> k_1,k_4> k_3, k_5> 2k_4-k_3\\
g_{k_1,k_1,k_3,k_3,k_3}=t^{k_3-k_1}\mathbb{L}^{-k_1}(\mathbb{L}-1)\sum_{k_2,k_4,k_5}g_{k_1,k_2,k_3,k_4,k_5}\\
g_{k_1,k_2,k_3,k_4,k_5}=g_{k_2,k_1,k_5,k_4,k_3}.\\
\end{cases}
\end{equation}

Moreover, $g_{k_1,k_2,k_3,k_4,k_5}=0$, if the collection $(k_1,k_2,k_3,k_4,k_5)$
does not satisfy the inequalities
$$k_2\ge k_1,k_4\ge k_3, k_5\ge 2k_4-k_3$$
or
$$k_1\ge k_2,k_4\ge k_5, k_3\ge 2k_4-k_5.$$

Similarly to the proof of Lemma 2, one can check that any coefficient could be expressed via
$g_{1,1,1,1,1}$ and the solution is unique.

Since by Theorem 2 $\int_{\mathcal{L}}t^{\mu}a^{v_x}b^{v_y}c^{v_x^2}d^{v_xv_y}f^{v_y^2}d\chi_{g}$ satisfies the equation (\ref{mu}), every solution of this equation is proportional to it.

Therefore $$I(t,a,b,c,d,f)=\sum G_{i,j}(t)a^{i}b^{j}c^{i^2}d^{ij}e^{j^2}.$$
Thus  (\ref{sys}) 
yields the partial differential equations (\ref{pde}).


Let us compute $g_{1,1,1,1,1}.$ If $v_x=v_y=1,$ then $\mu=0$,
so $$g_{1,1,1,1,1}=\chi_g\{v_x=v_y=1\}=(\mathbb{L}-1)^{2}\mathbb{L}^{-2}.$$

\section{ Examples} 

Recall that  $$G_{i,j}(t)=\int_{v_x=i,v_y=j}t^{\mu}d\chi_g.$$
One can check, that the system (\ref{sys}) gives a system of equations for  $G_{i,j}(t)$ of the form (\ref{recur}) with 
with $$\varepsilon_k=t^{k^2-k}\mathbb{L}^{-k}, C=(\mathbb{L}-1).$$
Therefore one has 
(from the proof of lemma 2) 
$$\begin{cases}
G_{i,j}(t)=t^{j^2-j}\mathbb{L}^{-j}G_{i-j,j},\ j<i\\
G_{i,i}(t)={(\mathbb{L}-1)t^{i^2-i}\mathbb{L}^{-i}\over 1-t^{i^2-i}\mathbb{L}^{1-i}}\sum_{j<i}G_{i,j}(t).\\
G_{1,1}(t)=(\mathbb{L}-1)^{2}\mathbb{L}^{-2}.$$
\end{cases}
$$
Then $G_{1,n}(t)=\mathbb{L}^{-1}G_{1,n-1},$ hence 
$$G_{1,n}(t)=(\mathbb{L}-1)^{2}\mathbb{L}^{-1-n}.$$
This corresponds to that fact that every arc with $v_x=1$ is smooth, so $\mu=0$.

Moreover, $G_{2,n}(t)=t^2\mathbb{L}^{-2}G_{2,n-2}(t)$ for $n>2$,
so  
$$G_{2,2n-1}=t^{2n-2}\mathbb{L}^{2-2n}G_{2,1}=
(\mathbb{L}-1)^{2}t^{2n-2}\mathbb{L}^{-1-2n}.$$
Every arc with $v_x=2$ and $v_y$ odd has the Milnor number equal to $v_y-1$, because the singularity is of $A_{v_y-1}$ type.

The result for the case when $v_y$ is even is more interesting.
$$G_{2,2n}(t)=t^{2n-2}\mathbb{L}^{2-2n}G_{2,2}=t^{2n}\mathbb{L}^{-2n}
{(\mathbb{L}-1)\over 1-t^{2}\mathbb{L}^{-1}}G_{1,2}(t)=$$
$${t^{2n}\mathbb{L}^{-2n-3}(\mathbb{L}-1)^3\over 1-t^{2}\mathbb{L}^{-1}}=
\sum_{k=0}^{\infty}t^{2(n+k)}\mathbb{L}^{-2n-3-k}(\mathbb{L}-1)^3.$$

Let us explain this answer. Making a change of variables one can obtain 
$x(t)=t^2$. Let $$y(t)=a_{2n}t^{2n}+\ldots+a_{2m+1}t^{2m+1}+\ldots,$$
where $a_{2m+1}$ is the first non-zero coefficient with odd index.
Then the equation of this curve is
$$F(x,y)=x^{2m+1}+(y-a_{2n}x^n-\ldots-a_{2m}x^{m})^2+\ldots=0,$$
so the Milnor number equals to $2m$. Consider $k=m-n$.
For fixed $x(t)$ the measure of series $y(t)$ with given $m$
is equal to $$(\mathbb{L}-1)^2\mathbb{L}^{-2n-(m-n+1}=(\mathbb{L}-1)^2\mathbb{L}^{-2n-k-1},$$ and the measure of series $x(t)$ with order 2 is equal to $(\mathbb{L}-1)\mathbb{L}^{-2}.$ Multiplying these expressions, we obtain the above formula. 

\begin{lemma}
Let a be the greatest common divisor of i and j. Then
$$G_{i,j}(t)=t^{(i-1)(j-1)-(a-1)^2}\mathbb{L}^{2a-i-j}G_{a,a}(t).$$
\end{lemma}

\begin{proof}
For $i=j$ the formula is tautological. Suppose that it is true for
$i-j$ and $j$ ($i>j$). Then 
$$G_{i,j}(t)=t^{j^2-j}\mathbb{L}^{-j}G_{i-j,j}=t^{(i-j-1)(j-1)-(a-1)^2+j(j-1)}\mathbb{L}^{2a-(i-j)-j-j}G_{a,a}(t).$$
Therefore the proposition follows from the Euclid algorithm.
\end{proof}

Similarly to the discussion above one can obtain
the following answers:

\vskip 1 cm
		\begin{tabular}{|l|l|l|l|l|}
			\hline
$a$ & 1 & 2 & 3 & 4  \\ \hline
& & & & \\
$G_{a,a}(t)$ & $(\mathbb{L}-1)\mathbb{L}^2$ & ${(\mathbb{L}-1)^3t^{2}\mathbb{L}^{-5}\over 1-t^{2}\mathbb{L}^{-1}}$ &
${(\mathbb{L}-1)^3t^{6}\mathbb{L}^{-7}(1+t^2\mathbb{L}^{-1})\over 
1-t^6\mathbb{L}^{-2}}$ & 
${(\mathbb{L}-1)^3t^{12}\mathbb{L}^{-9}(1-t^2\mathbb{L}^{-1}+t^{6}\mathbb{L}^{-1}-t^{8}\mathbb{L}^{-3})\over (1-t^{12}\mathbb{L}^{-3})(1-t^{2}\mathbb{L}^{-1})}$
\\ \hline
		\end{tabular}

\vskip 1 cm

This table together with Lemma 3  provides $G_{i,j}(t)$
with  $\mbox{\rm gcd}(i,j)\le 4.$

\begin{proposition}
Let $a=\mbox{\rm gcd}(i,j)$. $G_{ij}(t)$ is a power series of $t$, which coefficients
are Laurent polynomials in $\mathbb{L}$. If $a=1$, then
$$G_{i,j}(t)=(\mathbb{L}-1)^2t^{(i-1)(j-1)}\mathbb{L}^{-i-j},$$
and 
$$G_{i,j}(t)=(\mathbb{L}-1)^3t^{(i-1)(j-1)+a-1}\mathbb{L}^{-i-j-1}+\mbox{\rm terms of higher degree in $t$}.$$
\end{proposition}
\begin{proof}
The first statement can be easily checked by induction. The case $a=1$ follows from Lemma 3. Let us prove the formula for the case $a>1$.
One has 
$$G_{a,a}(t)={(\mathbb{L}-1)t^{a^2-a}\mathbb{L}^{-a}\over 1-t^{a^2-a}\mathbb{L}^{1-a}}(G_{a,1}(t)+\mbox{\rm O}(t))=$$
$$(\mathbb{L}-1)t^{a^2-a}\mathbb{L}^{-a}(1+\mbox{\rm O}(t))((\mathbb{L}-1)^2\mathbb{L}^{-a-1}+\mbox{\rm O}(t))=
(\mathbb{L}-1)^{3}t^{a^2-a}\mathbb{L}^{-2a-1}+\mbox{\rm O}(t^{a^2-a+1}) .$$
Now the statement follows from Lemma 3.
\end{proof}

\section{Functional equation for the intersection number}

Let 
$${\scriptstyle
J(t,a,b,c,d,p,q,r,s)=\int_{\mathcal{L}^{(1)}\times\mathcal{L}^{(2)}}t^{\gamma_1\circ\gamma_2}a^{v_x^{(1)}v_x^{(2)}}
b^{v_x^{(1)}v_y^{(2)}}c^{v_y^{(1)}v_x^{(2)}}d^{v_y^{(1)}v_y^{(2)}}p^{v_x^{(1)}}
q^{v_y^{(1)}}r^{v_x^{(2)}}s^{v_y^{(2)}}d\chi_g,}$$

$${\scriptstyle
B(t,a,b,c,d,p,q,r,s)=\int_{\{v_y^{(1)}>v_x^{(1)},v_y^{(2)}>v_x^{(2)}\}}t^{\gamma_1\circ\gamma_2}a^{v_x^{(1)}v_x^{(2)}}
b^{v_x^{(1)}v_y^{(2)}}c^{v_y^{(1)}v_x^{(2)}}d^{v_y^{(1)}v_y^{(2)}}p^{v_x^{(1)}}
q^{v_y^{(1)}}r^{v_x^{(2)}}s^{v_y^{(2)}}d\chi_g,}$$

Note that $$\sigma(\gamma_1)\circ\sigma(\gamma_2)=\gamma_1\circ\gamma_2+v_1v_2,$$
hence if after blowing up the origin arcs intersect the exceptional
divisor at different points, their intersection number 
equals to the product of multiplicities.

Let us decompose $\mathcal{L}^{(1)}\times\mathcal{L}^{(2)}$ into components with respect to the inequalities between
$v_y^{(1)}$ and $v_x^{(1)}$, and also $v_y^{(2)}$ and $v_x^{(2)}.$

1)$v_y^{(1)}>v_x^{(1)},v_y^{(2)}>v_x^{(2)}.$ By definition, the integral equals to
$$B(t,a,b,c,d,p,q,r,s).$$

2)$v_y^{(1)}>v_x^{(1)},v_y^{(2)}<v_x^{(2)}.$ After blow-up arcs intersect the exceptional divisor at different
points, and therefore the integral equals to

$$\int_{\{v_y^{(1)}>v_x^{(1)},v_y^{(2)}<v_x^{(2)}\}}t^{v_x^{(1)}v_y^{(2)}}a^{v_x^{(1)}v_x^{(2)}}
b^{v_x^{(1)}v_y^{(2)}}c^{v_y^{(1)}v_x^{(2)}}d^{v_y^{(1)}v_y^{(2)}}p^{v_x^{(1)}}
q^{v_y^{(1)}}r^{v_x^{(2)}}s^{v_y^{(2)}}d\chi_g=$$
$$(v_x^{(1)}=i,v_y^{(1)}=j,v_y^{(2)}=k,v_x^{(2)}=l)=$$
$$\sum_{i<j,k<l}a^{il}(bt)^{ik}c^{jl}d^{jk}p^{i}q^{j}r^{l}s^{k}(\mathbb{L}-1)^{4}\mathbb{L}^{-i-j-k-l}=$$
$$(\mathbb{L}-1)^4\Phi(bt,a,d,c;\mathbb{L}^{-1}p,\mathbb{L}^{-1}q,\mathbb{L}^{-1}s,\mathbb{L}^{-1}r),$$ 
where  $\Phi$ is defined by the equation
$$\Phi(\alpha,\beta,\gamma,\delta,\pi,\kappa,\rho,\sigma)=\sum_{i<j,k<l}\alpha^{ik}\beta^{il}\gamma^{jk}\delta^{jl}\pi^{i}\kappa^{j}\rho^{k}\sigma^{l}.$$

3)$v_y^{(1)}>v_x^{(1)},v_y^{(2)}=v_x^{(2)}.$ After blowing up arcs intersect the exceptional divisor at different
points. Therefore the integral equals to
$$\int_{\{v_y^{(1)}>v_x^{(1)},v_y^{(2)}=v_x^{(2)}\}}t^{v_x^{(1)}v_x^{(2)}}a^{v_x^{(1)}v_x^{(2)}}
b^{v_x^{(1)}v_x^{(2)}}c^{v_y^{(1)}v_x^{(2)}}d^{v_y^{(1)}v_x^{(2)}}p^{v_x^{(1)}}
q^{v_y^{(1)}}r^{v_x^{(2)}}s^{v_x^{(2)}}d\chi_g=$$
$$\sum_{i<j,k}(tab)^{ik}(cd)^{jk}p^{i}q^{j}(rs)^{k}(\mathbb{L}-1)^4\mathbb{L}^{-i-j-2k}=(\mathbb{L}-1)^4\Psi(tab,cd,p\mathbb{L}^{-1},q\mathbb{L}^{-1},rs\mathbb{L}^{-2}),$$
where $\Psi$ is defined by the equation
$$\Psi(\alpha,\beta,\pi,\kappa,\rho)=\sum_{i<j,k}\alpha^{ik}\beta^{jk}\pi^{i}\kappa^{j}\rho^{k}.$$

4)$v_y^{(1)}<v_x^{(1)},v_y^{(2)}>v_x^{(2)}.$
$$\int_{\{v_y^{(1)}<v_x^{(1)},v_y^{(2)}>v_x^{(2)}\}}t^{v_{y}^{(1)}v_x^{(2)}}a^{v_x^{(1)}v_x^{(2)}}
b^{v_x^{(1)}v_y^{(2)}}c^{v_y^{(1)}v_x^{(2)}}d^{v_y^{(1)}v_y^{(2)}}p^{v_x^{(1)}}
q^{v_y^{(1)}}r^{v_x^{(2)}}s^{v_y^{(2)}}d\chi_g=$$
$$\sum_{i<j,k<l}a^{jk}b^{jl}(ct)^{ik}d^{il}p^{j}q^{i}r^{k}s^{l}(\mathbb{L}-1)^{4}\mathbb{L}^{-i-j-k-l}=$$
$$(\mathbb{L}-1)^4\Phi(ct,d,a,b,q\mathbb{L}^{-1},p\mathbb{L}^{-1},r\mathbb{L}^{-1},s\mathbb{L}^{-1}).$$

5)$v_y^{(1)}<v_x^{(1)},v_y^{(2)}<v_x^{(2)}.$ From the symmetry it is clear that the integral
equals to
$$\int_{\{v_y^{(1)}>v_x^{(1)},v_y^{(2)}>v_x^{(2)}\}}t^{\gamma_1\circ\gamma_2}a^{v_y^{(1)}v_y^{(2)}}
b^{v_y^{(1)}v_x^{(2)}}c^{v_x^{(1)}v_y^{(2)}}d^{v_x^{(1)}v_x^{(2)}}p^{v_y^{(1)}}
q^{v_x^{(1)}}r^{v_y^{(2)}}s^{v_x^{(2)}}d\chi_g=B(t,d,b,c,a,q,p,s,r).$$

6)$v_y^{(1)}<v_x^{(1)},v_y^{(2)}=v_x^{(2)}.$
$$\int_{\{v_y^{(1)}<v_x^{(1)},v_y^{(2)}=v_x^{(2)}\}}t^{v_y^{(1)}v_x^{(2)}}a^{v_x^{(1)}v_x^{(2)}}
b^{v_x^{(1)}v_x^{(2)}}c^{v_y^{(1)}v_x^{(2)}}d^{v_y^{(1)}v_x^{(2)}}p^{v_x^{(1)}}
q^{v_y^{(1)}}r^{v_x^{(2)}}s^{v_x^{(2)}}d\chi_g=$$ 
$$\sum_{i<j,k}(tcd)^{ik}(ab)^{jk}p^{j}q^{i}(rs)^{k}(\mathbb{L}-1)^4\mathbb{L}^{-i-j-2k}=
(\mathbb{L}-1)^4\Psi(tcd,ab,q\mathbb{L}^{-1},p\mathbb{L}^{-1},rs\mathbb{L}^{-2}).$$

7)$v_y^{(1)}=v_x^{(1)},v_y^{(2)}>v_x^{(2)}.$
$$\int_{\{v_y^{(1)}=v_x^{(1)},v_y^{(2)}>v_x^{(2)}\}}t^{v_x^{(1)}v_x^{(2)}}a^{v_x^{(1)}v_x^{(2)}}
b^{v_x^{(1)}v_y^{(2)}}c^{v_x^{(1)}v_x^{(2)}}d^{v_x^{(1)}v_y^{(2)}}p^{v_x^{(1)}}
q^{v_x^{(1)}}r^{v_x^{(2)}}s^{v_y^{(2)}}d\chi_g=$$ 
$$\sum_{i<j,k}(tbd)^{ik}(ac)^{jk}(pq)^{k}r^{i}s^{j}(\mathbb{L}-1)^4\mathbb{L}^{-i-j-2k}=
(\mathbb{L}-1)^4\Psi(tac,bd,r\mathbb{L}^{-1},s\mathbb{L}^{-1},pq\mathbb{L}^{-2}).$$

8)$v_y^{(1)}=v_x^{(1)},v_y^{(2)}<v_x^{(2)}.$
$$\int_{v_y^{(1)}=v_x^{(1)},v_y^{(2)}<v_x^{(2)}}t^{v_x^{(1)}v_x^{(2)}}a^{v_x^{(1)}v_x^{(2)}}
b^{v_x^{(1)}v_y^{(2)}}c^{v_x^{(1)}v_x^{(2)}}d^{v_x^{(1)}v_y{(2)}}p^{v_x^{(1)}}
q^{v_x^{(1)}}r^{v_x^{(2)}}s^{v_y^{(2)}}d\chi_g=$$ 
$$\sum_{i<j,k}(tbd)^{ik}(ac)^{jk}r^{j}s^{i}(pq)^{k}(\mathbb{L}-1)^4\mathbb{L}^{-i-j-2k}=
(\mathbb{L}-1)^4\Psi(tbd,ac,s\mathbb{L}^{-1},r\mathbb{L}^{-1},pq\mathbb{L}^{-2}).$$

9)$v_y^{(1)}=v_x^{(1)},v_y^{(2)}=v_x^{(2)}$. In this case
$$\begin{cases}
y^{(1)}=\lambda_{1}x^{(1)}+\widetilde{y}^{(1)}, v_{\widetilde{y}}^{(1)}>v_x^{(1)},\\ 
y^{(2)}=\lambda_{2}x^{(2)}+\widetilde{y}^{(2)}, v_{\widetilde{y}}^{(2)}>v_x^{(2)},\\
\end{cases}$$
$\lambda_1\neq 0,\lambda_2\neq 0.$

a)$\lambda_1\neq\lambda_2.$ In this case
$$[\{\lambda_1,\lambda_2\in\mathbb{C}^{*}|\lambda_1\neq\lambda_2\}]=
[\mathbb{C}^{*}]^2-[\{\lambda_1,\lambda_2\in\mathbb{C}^{*}|\lambda_1=\lambda_2\}]=
(\mathbb{L}-1)^2-(\mathbb{L}-1)=(\mathbb{L}-1)(\mathbb{L}-2),$$
and for fixed $(\lambda_1,\lambda_2)$  the integral equals to
$$\int_{\{v_{\widetilde{y}}^{(1)}>v_x^{(1)},v_{\widetilde{y}}^{(2)}>v_x^{(2)}\}}t^{v_x^{(1)}v_x^{(2)}}a^{v_x^{(1)}v_x^{(2)}}
b^{v_x^{(1)}v_x^{(2)}}c^{v_x^{(1)}v_x^{(2)}}d^{v_x^{(1)}v_x^{(2)}}p^{v_x^{(1)}}
q^{v_x^{(1)}}r^{v_x^{(2)}}s^{v_x^{(2)}}d\chi_g=$$
$$\sum_{i<j,k<l}(tabcd)^{ik}(pq)^{i}(rs)^{k}(\mathbb{L}-1)^{4}\mathbb{L}^{-i-j-k-l}=$$
$$(\mathbb{L}-1)^4\Phi(tabcd,1,1,1;\mathbb{L}^{-1}pq,\mathbb{L}^{-1},\mathbb{L}^{-1}rs,\mathbb{L}^{-1}).$$
Therefore the contribution of the stratum
 $\{v_x^{(1)}=v_y^{(1)},v_x^{(2)}=v_y^{(2)},\lambda_1\neq\lambda_2\}$
equals to
$$(\mathbb{L}-1)^5(\mathbb{L}-2)\Phi(tabcd,1,1,1;\mathbb{L}^{-1}pq,\mathbb{L}^{-1},\mathbb{L}^{-1}rs,\mathbb{L}^{-1}).$$

b)$\lambda_1=\lambda_2=\lambda.$ Under the affine change of variables
$(x,y)\mapsto (x,y-\lambda x)$ intersection number does not change,
so for fixed $\lambda$ the integral equals to
$$\int_{\{v_{\widetilde{y}}^{(1)}>v_x^{(1)},v_{\widetilde{y}}^{(2)}>v_x^{(2)}\}}t^{\gamma_1\circ\gamma_2}a^{v_x^{(1)}v_x^{(2)}}
b^{v_x^{(1)}v_x^{(2)}}c^{v_x^{(1)}v_x^{(2)}}d^{v_x^{(1)}v_x^{(2)}}p^{v_x^{(1)}}
q^{v_x^{(1)}}r^{v_x^{(2)}}s^{v_x^{(2)}}d\chi_g=$$
$$B(t,abcd,1,1,1,pq,1,rs,1),$$
hence the total contribution of this stratum equals to
$$(\mathbb{L}-1)B(t,abcd,1,1,1,pq,1,rs,1).$$

So we have:

$${\scriptstyle
J(t,a,b,c,d,p,q,r,s)=B(t,a,b,c,d,p,q,r,s)+
(\mathbb{L}-1)^4\Phi(bt,a,d,c;\mathbb{L}^{-1}p,\mathbb{L}^{-1}q,\mathbb{L}^{-1}s,\mathbb{L}^{-1}r)+}$$
$${\scriptstyle
(\mathbb{L}-1)^4\Psi(tab,cd,p\mathbb{L}^{-1},q\mathbb{L}^{-1},rs\mathbb{L}^{-2})+
(\mathbb{L}-1)^4\Phi(ct,d,a,b,q\mathbb{L}^{-1},p\mathbb{L}^{-1},r\mathbb{L}^{-1},s\mathbb{L}^{-1})+}$$
\begin{equation}
\label{inters}
{\scriptstyle
B(t,d,c,b,a,q,p,s,r)+(\mathbb{L}-1)^4\Psi(tcd,ab,q\mathbb{L}^{-1},p\mathbb{L}^{-1},rs\mathbb{L}^{-2})+}
\end{equation}
$${\scriptstyle
(\mathbb{L}-1)^4\Psi(tbd,ac,s\mathbb{L}^{-1},r\mathbb{L}^{-1},pq\mathbb{L}^{-2})+
(\mathbb{L}-1)^4\Psi(tac,bd,r\mathbb{L}^{-1},s\mathbb{L}^{-1},pq\mathbb{L}^{-2})+}$$
$${\scriptstyle
(\mathbb{L}-1)^5(\mathbb{L}-2)\Phi(tabcd,1,1,1;\mathbb{L}^{-1}pq,\mathbb{L}^{-1},\mathbb{L}^{-1}rs,\mathbb{L}^{-1})+
(\mathbb{L}-1)B(t,abcd,1,1,1,pq,1,rs,1).}$$

Consider a blow-up. If $v_y>v_x$, then $y(t)=x(t)\theta(t)$.
So $$\sigma^{*}(\gamma_1\circ\gamma_2)=\gamma_1\circ\gamma_2+v_x^{(1)}v_x^{(2)},\sigma^{*}v_x^{(1)}=v_x^{(1)},$$
$$\sigma^{*}v_x^{(2)}=v_x^{(2)},\sigma^{*}v_y^{(1)}=v_x^{(1)}+v_{\theta}^{(1)},\sigma^{*}v_y^{(2)}=v_x^{(2)}+v_{\theta}^{(2)}.$$
Using the Denef-Loeser change of variables formula we have
$$B(t,a,b,c,d,p,q,r,s)=\int_{\mathcal{L}^{(1)}\times\mathcal{L}^{(2)}}t^{\gamma_1\circ\gamma_2+v_x^{(1)}v_x^{(2)}}a^{v_x^{(1)}v_x^{(2)}}
b^{v_x^{(1)}v_x^{(2)}+v_x^{(1)}v_{\theta}^{(2)}}c^{v_x^{(1)}v_x^{(2)}+v_{\theta}^{(1)}v_x^{(2)}}\times$$
$$\times d^{v_x^{(1)}v_x^{(2)}+v_{\theta}^{(1)}v_x^{(2)}+v_x^{(1)}v_{\theta}^{(2)}+v_{\theta}^{(1)}v_{\theta}^{(2)}}p^{v_x^{(1)}}
q^{v_x^{(1)}+v_{\theta}^{(1)}}r^{v_x^{(2)}}s^{v_x^{(2)}+v_{\theta}^{(2)}}\mathbb{L}^{-v_x^{(1)}-v_x^{(2)}}d\chi_g=$$
$$J(t,tabcd,bd,cd,d,pq\mathbb{L}^{-1},q,rs\mathbb{L},s).$$

Substituting these expression for $B(\cdot)$ into (\ref{inters}), we get the following statement.

\begin{lemma}
The function  $J$ satisfies the functional equation: 
$${\scriptstyle
J(t,a,b,c,d,p,q,r,s)=J(t,tabcd,bd,cd,d,pq\mathbb{L}^{-1},q,rs\mathbb{L}^{-1},s)+}$$
$${\scriptstyle
(\mathbb{L}-1)^4\Phi(bt,a,d,c,\mathbb{L}^{-1}p,\mathbb{L}^{-1}q,\mathbb{L}^{-1}s,\mathbb{L}^{-1}r)+
(\mathbb{L}-1)^4\Psi(tab,cd,p\mathbb{L}^{-1},q\mathbb{L}^{-1},rs\mathbb{L}^{-2})+}$$
$${\scriptstyle
(\mathbb{L}-1)^4\Phi(ct,d,a,b,q\mathbb{L}^{-1},p\mathbb{L}^{-1},r\mathbb{L}^{-1},s\mathbb{L}^{-1})+
J(t,tabcd,ac,ab,a,pq\mathbb{L}^{-1},p,rs\mathbb{L}^{-1},r)}+$$
$${\scriptstyle
(\mathbb{L}-1)^4\Psi(tcd,ab,q\mathbb{L}^{-1},p\mathbb{L}^{-1},rs\mathbb{L}^{-2})+
(\mathbb{L}-1)^4\Psi(tbd,ac,s\mathbb{L}^{-1},r\mathbb{L}^{-1},pq\mathbb{L}^{-2})+}$$
$${\scriptstyle
(\mathbb{L}-1)^4\Psi(tac,bd,r\mathbb{L}^{-1},s\mathbb{L}^{-1},pq\mathbb{L}^{-2})+
(\mathbb{L}-1)^5(\mathbb{L}-2)\Phi(tabcd,1,1,1;\mathbb{L}^{-1}pq,\mathbb{L}^{-1},\mathbb{L}^{-1}rs,\mathbb{L}^{-1})+}$$
$${\scriptstyle
(\mathbb{L}-1)J(t,tabcd,1,1,1,pq\mathbb{L}^{-1},1,rs\mathbb{L}^{-1},1).}$$

\end{lemma}

\section{ Power structures}

The notion of the power structure over a (semi)ring was introduced
by S. Gusein-Zade, I. Luengo and A. Melle-Hernandez in \cite{powers}.

\begin{definition}
A power structure on the ring $R$
is a map $$(1+tR[[t]])\times R\rightarrow 1+tR[[t]]:(A(t),m)\mapsto (A(t))^{m},$$
satisfying the following properties:

1.$(A(t))^{0}=1,$

2.$(A(t))^{1}=A(t),$

3.$((A(t)\cdot B(t))^{m}=((A(t))^{m}\cdot((B(t))^{m},$

4.$(A(t))^{m+n}=(A(t))^{m}\cdot(A(t))^{n},$

5.$(A(t))^{mn}=((A(t))^{n})^{m},$

6.$(1+t)^{m}=1+mt+$ terms with higher degree,

7.$(A(t^{k}))^{m}=((A(t))^{m})|_{t\rightarrow t^{k}}.$
\end{definition}

A power structure is called finitely determined
if for every $N>0$ there exists such $M>0$ that the
$N$-jet of the series $(A(t))^{m}$ is uniquely determined by the
$M$-jet of the series $A(t)$.

To fix a finitely determined power structure
it is sufficient to define the series $(1-t)^{-m}$
for every $m\in R$ such that $(1-t)^{-m-n}=(1-t)^{-m}\cdot(1-t)^{-n}$.
Over the Grothendieck ring of varieties there is 
the power structure, such that 
$$(1-t)^{-[X]}=1+[S^{1}X]t+[S^{2}X]t^{2}+\ldots,$$
where $S^{k}X=X^{k}/S_{k}$ denotes the $k$-th symmetric power of $X$.
For example, for $j\ge 0$ 
$$(1-t)^{-\mathbb{L}^j}=\sum_{k=0}^{\infty}t^k\mathbb{L}^{kj}=(1-t\mathbb{L}^j)^{-1}.$$

For $X\in K_{0}(Var_{\mathbb{C}})$, $k>0$, let
$$(1-t)^{-\mathbb{L}^{-k}X}=(1-u)^{-X}|_{u=\mathbb{L}^{-k}t}.$$

The following statement defines the corresponding power structure over the ring $\mathcal{M}$

\begin{lemma}
The map $K_0(Var_{\mathbb{C}})[\mathbb{L}^{-1}]
\rightarrow 1+tK_0(Var_{\mathbb{C}})[\mathbb{L}^{-1}][[t]],$
$$Z\mapsto (1-t)^{-Z}$$ is well defined. It transforms the sum
into the product and is continuous with respect to the filtration $F_k$.
\end{lemma}




Let us construct a measure on the symmetric power $S^{k}\mathcal{L}$.
For a cylindric set $A=\pi_{n}^{-1}(A_n)$
let $\mu(S^{k}A)=\mathbb{L}^{-2nk}[S^{k}A_{n}]$.
This construction corresponds to the power structure
over the ring $\mathcal{M}$,
so $$\sum_{k}\mu(S^{k}A)t^{k}=(1-t)^{-\mu(A)}.$$

Let $B_{i}$ be a collection of non-intersecting cylindric subsets of
$\mathcal{L}$,let $k_i$ be nonnegative integers with 
 $\sum k_i=n$. Then for the natural embedding
$$S^{k_1}B_1\times S^{k_2}B_2\times\ldots\rightarrow S^{n}\mathcal{L}$$ 
let $\mu(\prod_{i}S^{k_i}B_{i})=\prod_{i}\mu(S_{k_i}B_{i}).$
Consider the algebra of sets generated by such products
of symmetric powers of cylindric sets.
Continuation of 
$\mu$ is a well defined additive measure on this algebra.

\begin{lemma}
Let $f$ be a simple function on $\mathcal{L}$.
Define a function $F$ on $\sqcup_{k}S^{k}\mathcal{L}$
by the formula $F(\gamma_1,\ldots,\gamma_k)=\prod_{i}f(\gamma_{i}).$ 
Then
$$\int_{\sqcup_{k}S^{k}\mathcal{L}}Fd\chi_g=\int_{\mathcal{L}}(1-f)^{-d\chi_g}.$$
Here $d\chi_g$ is in the exponent to emphasize that $1-f$ is considered as an element of an abelian group with respect to multiplication. 
\end{lemma}

\begin{proof}
Consider $\mathcal{L}=\sqcup_{j}B_{j}, f|_{B_j}=f_{j}.$
Then
$$S^{k}\mathcal{L}=\sqcup_{k_1+k_2+\ldots=k}S^{k_1}B_{1}\times S^{k_2}B_{2}\times\ldots,$$
and
$$F|_{S^{k_1}B_{1}\times S^{k_2}B_{2}\times\ldots}=f_{1}^{k_1}\cdot f_{2}^{k_2}\cdot\ldots.$$
Therefore
$$\int_{\sqcup_{k}S^{k}\mathcal{L}}=\sum_{k}\sum_{\sum k_i=k}\mu(S^{k_1}B_{1}\times S^{k_2}B_2\times \ldots)f_{1}^{k_1}f_{2}^{k_2}\ldots=$$
$$\prod_{j}(1+\mu(S^{1}B_{j})f_{j}+\mu(S^{2}B_{j}){f_j}^2+\ldots)=\prod_{j}(1-f_j)^{-\mu(B_{j})}=\int_{\mathcal{L}}(1-f)^{-d\chi_g}.$$
\end{proof}

\section{An equation for the intersection numbers on $$\mathcal{L}\times\sqcup_{k}S^{k}\mathcal{L}$$}

Consider the generating function:
$${\scriptstyle
I(t,a,b,c,d,p,q,r,s,u)=\int_{\mathcal{L}\times\sqcup_{k}S^{k}\mathcal{L}}t^{\gamma_1\circ\gamma_2}a^{v_x^{(1)}v_x^{(2)}}
b^{v_x^{(1)}v_y^{(2)}}c^{v_y^{(1)}v_x^{(2)}}d^{v_y^{(1)}v_y^{(2)}}p^{v_x^{(1)}}
q^{v_y^{(1)}}r^{v_x^{(2)}}s^{v_y^{(2)}}u^{k}d\chi_g.}$$
In this section we will obtain a functional equation for the function $I$.

Let $$\prod_{k\ge l}(1-x^{k}y^{l}u)^{-(\mathbb{L}-1)^2}=\sum_{k_1,k_2}\varepsilon_{k_1,k_2}(u)x^{k_1}y^{k_2}$$
and
$$\prod_{k<l}(1-(xy)^{k}z^{l}u)^{-(\mathbb{L}-1)^2}\prod_{k>l}(1-x^{k}(yz)^{l}u)^{-(\mathbb{L}-1)^2}
\prod_{k<l}(1-(xyz)^{k}\mathbb{L}^{-l}u)^{-(\mathbb{L}-2)(\mathbb{L}-1)^2}=$$
$$\sum_{k_1,k_2,k_3}\alpha_{k_1,k_2,k_3}(u)x^{k_1}y^{k_2}z^{k_3}.$$
Consider

$$J_{\gamma_1}(t,a,b,c,d,r,s,u)=\int_{\sqcup_{k}S^{k}\mathcal{L}}t^{\gamma_1\circ\gamma_2}a^{v_x^{(1)}v_x^{(2)}}
b^{v_x^{(1)}v_y^{(2)}}c^{v_y^{(1)}v_x^{(2)}}d^{v_y^{(1)}v_y^{(2)}}
r^{v_x^{(2)}}s^{v_y^{(2)}}u^{k}d\chi_g=$$
$$\int_{\mathcal{L}^{(2)}}
(1-t^{\gamma_1\circ\gamma_2}a^{v_x^{(1)}v_x^{(2)}}
b^{v_x^{(1)}v_y^{(2)}}c^{v_y^{(1)}v_x^{(2)}}d^{v_y^{(1)}v_y^{(2)}}
r^{v_x^{(2)}}s^{v_y^{(2)}}u)^{-d\chi_g^{(2)}},$$
then 
$$I(t,a,b,c,d,p,q,r,s,u)=\int_{\mathcal{L}^{(1)}}p^{v_x^{(1)}}q^{v_y^{(1)}}J_{\gamma_1}(t,a,b,c,d,r,s,u)d\chi_g^{(1)}.$$
If $v_y^{(1)}>v_x^{(1)},$ then
$$J_{\gamma_1}(t,a,b,c,d,r,s,u)=
\int_{\{v_y^{(2)}>v_x^{(2)}\}}(1-t^{\gamma_1\circ\gamma_2}a^{v_x^{(1)}v_x^{(2)}}
b^{v_x^{(1)}v_y^{(2)}}c^{v_y^{(1)}v_x^{(2)}}d^{v_y^{(1)}v_y^{(2)}}
r^{v_x^{(2)}}s^{v_y^{(2)}}u)^{-d\chi_g^{(2)}}\times$$
$$\int_{\{v_y^{(2)}\le v_x^{(2)}\}}(1-t^{v_x^{(1)}v_y^{(2)}}a^{v_x^{(1)}v_x^{(2)}}
b^{v_x^{(1)}v_y^{(2)}}c^{v_y^{(1)}v_x^{(2)}}d^{v_y^{(1)}v_y^{(2)}}
r^{v_x^{(2)}}s^{v_y^{(2)}}u)^{-d\chi_g^{(2)}}.$$
Let us blow-up the origin.
By the Denef-Loeser formula
$d\chi_g^{(2)}\mapsto\mathbb{L}^{-v_x^{(2)}}d\chi_g^{(2)}$,
so the first integral is equal to 
$$\int_{\mathcal{L}^{(2)}}(1-t^{\sigma^{-1}(\gamma_1)\circ\gamma_2}(tabcd)^{\widetilde{v_x}^{(1)}v_x^{(2)}}
(bd)^{\widetilde{v_x}^{(1)}v_y^{(2)}}(cd)^{\widetilde{v_y}^{(1)}v_x^{(2)}}d^{\widetilde{v_y}^{(1)}v_y^{(2)}} 
(rs)^{v_x^{(2)}}s^{v_y^{(2)}}u)^{-\mathbb{L}^{-v_x^{(2)}}d\chi_g^{(2)}}=$$
$$J_{\sigma^{-1}(\gamma_1)}(t,tabcd,bd,cd,d,rs\mathbb{L}^{-1},s,u).$$
The second integral equals to
$${\scriptstyle
\prod_{k\le l}(1-a^{lv_x^{(1)}}(bt)^{kv_x^{(1)}}c^{lv_y^{(1)}}d^{kv_y^{(1)}}r^{l}s^{k}u)^{-\mathbb{L}^{-k-l}(\mathbb{L}-1)^2}=
\prod_{k\le l}(1-a^{lv_x^{(1)}}(bt)^{kv_x^{(1)}}c^{lv_y^{(1)}}d^{kv_y^{(1)}}(\mathbb{L}^{-1}r)^{l}(\mathbb{L}^{-1}s)^{k}u)^{-(\mathbb{L}-1)^2}=}$$
$$=\sum_{k_1,k_2}\varepsilon_{k_1,k_2}(u)(\mathbb{L}^{-1}r)^{k_1}(\mathbb{L}^{-1}s)^{k_2}a^{k_1v_x^{(1)}}(bt)^{k_2v_x^{(1)}}c^{k_1v_y^{(1)}}d^{k_2v_y^{(1)}}.$$
Hence 
$$J_{\gamma_1}(t,a,b,c,d,r,s,u)=J_{\sigma^{-1}(\gamma_1)}(t,tabcd,bd,cd,d,rs\mathbb{L}^{-1},s,u)\times$$
$$(\sum_{k_1,k_2}\varepsilon_{k_1,k_2}(u)(\mathbb{L}^{-1}r)^{k_1}(\mathbb{L}^{-1}s)^{k_2}a^{k_1v_x^{(1)}}(bt)^{k_2v_x^{(1)}}c^{k_1v_y^{(1)}}d^{k_2v_y^{(1)}}).$$
Therefore we get that 
$${\scriptstyle
E(t,a,b,c,d,p,q,r,s,u)=\int_{\{v_y^{(1)}>v_x^{(1)}\}}t^{\gamma_1\circ\gamma_2}a^{v_x^{(1)}v_x^{(2)}}
b^{v_x^{(1)}v_y^{(2)}}c^{v_y^{(1)}v_x^{(2)}}d^{v_y^{(1)}v_y^{(2)}}p^{v_x^{(1)}}
q^{v_y^{(1)}}r^{v_x^{(2)}}s^{v_y^{(2)}}u^{k}d\chi_g=}$$
$${\scriptstyle
\int_{\{v_y^{(1)}>v_x^{(1)}\}}J_{\gamma_1}(t,a,b,c,d,r,s,u)p^{v_x^{(1)}}q^{v_y^{(1)}}d\chi_g^{(1)}=
\sum_{k_1,k_2}\varepsilon_{k_1,k_2}(u)(\mathbb{L}^{-1}r)^{k_1}(\mathbb{L}^{-1}s)^{k_2}\times}$$
$${\scriptstyle
\int_{\{v_y^{(1)}>v_x^{(1)}\}}J_{\sigma^{-1}(\gamma_1)}(t,tabcd,bd,cd,d,rs\mathbb{L}^{-1},s,u)a^{k_1v_x^{(1)}}(bt)^{k_2v_x^{(1)}}c^{k_1v_y^{(1)}}d^{k_2v_y^{(1)}}
p^{v_x^{(1)}}q^{v_y^{(2)}}d\chi_g^{(1)}=}$$
$${\scriptstyle
\sum_{k_1,k_2}\varepsilon_{k_1,k_2}(u)(\mathbb{L}^{-1}r)^{k_1}(\mathbb{L}^{-1}s)^{k_2}
\int_{\{v_y^{(1)}>v_x^{(1)}\}}J_{\sigma^{-1}(\gamma_1)}(t,tabcd,bd,cd,d,rs\mathbb{L}^{-1},s,u)\times}$$
$${\scriptstyle
 (a^{k_1}(bt)^{k_2}p)^{v_x^{(1)}}
(c^{k_1}d^{k_2}q)^{v_y^{(1)}}d\chi_g^{(1)}=
\sum_{k_1,k_2}\varepsilon_{k_1,k_2}(u)(\mathbb{L}^{-1}r)^{k_1}(\mathbb{L}^{-1}s)^{k_2}\times}$$
$${\scriptstyle
\int_{\mathcal{L}^{(1)}}J_{\gamma_1}(t,tabcd,bd,cd,d,rs\mathbb{L}^{-1},s,u)((ac)^{k_1}(bdt)^{k_2}pq\mathbb{L}^{-1})^{v_x^{(1)}}
(c^{k_1}d^{k_2}q)^{v_y^{(1)}}d\chi_g^{(1)}=}$$
$${\scriptstyle
\sum_{k_1,k_2}\varepsilon_{k_1,k_2}(u)(\mathbb{L}^{-1}r)^{k_1}(\mathbb{L}^{-1}s)^{k_2}
I(t,tabcd,bd,cd,d,(ac)^{k_1}(bdt)^{k_2}pq\mathbb{L}^{-1},c^{k_1}d^{k_2}q,rs\mathbb{L}^{-1},s,u).}$$
It is clear that the same integral over $\{v_x^{(1)}<v_y^{(1)}\}$
equals to $E(t,d,c,b,a,q,p,s,r,u)$.

Let us compute the integral over  $\{v_x^{(1)}=v_y^{(1)}\}.$ Let
$y^{(1)}=\lambda_1 x^{(1)}+\widetilde{y}^{(1)}.$
Then 
$${\scriptstyle
J_{\gamma_1}(t,a,b,c,d,r,s,u)=
\int_{\{v_x^{(2)}<v_y^{(2)}\}}(1-t^{\gamma_1\circ\gamma_2}a^{v_x^{(1)}v_x^{(2)}}
b^{v_x^{(1)}v_y^{(2)}}c^{v_y^{(1)}v_x^{(2)}}d^{v_y^{(1)}v_y^{(2)}}
r^{v_x^{(2)}}s^{v_y^{(2)}}u)^{-d\chi_g^{(2)}}\times}$$
$$\int_{\{v_x^{(2)}>v_y^{(2)}\}}(1-t^{\gamma_1\circ\gamma_2}a^{v_x^{(1)}v_x^{(2)}}
b^{v_x^{(1)}v_y^{(2)}}c^{v_y^{(1)}v_x^{(2)}}d^{v_y^{(1)}v_y^{(2)}}
r^{v_x^{(2)}}s^{v_y^{(2)}}u)^{-d\chi_g^{(2)}}\times$$
$$\int_{\{v_x^{(2)}=v_y^{(2)}\}}(1-t^{\gamma_1\circ\gamma_2}a^{v_x^{(1)}v_x^{(2)}}
b^{v_x^{(1)}v_y^{(2)}}c^{v_y^{(1)}v_x^{(2)}}d^{v_y^{(1)}v_y^{(2)}}
r^{v_x^{(2)}}s^{v_y^{(2)}}u)^{-d\chi_g^{(2)}}=$$
$${\scriptstyle
\prod_{k<l}(1-(tac)^{kv_x^{(1)}}(bd)^{lv_x^{(1)}}r^ks^{l}u)^{-\mathbb{L}^{-k-l}(\mathbb{L}-1)^2}\cdot
\prod_{k>l}(1-(ac)^{kv_x^{(1)}}(tbd)^{lv_x^{(1)}}r^ks^{l}u)^{-\mathbb{L}^{-k-l}(\mathbb{L}-1)^2}\times}$$
$$\int_{\{v_x^{(2)}<v_{\widetilde{y}}^{(2)}\}}(1-t^{\gamma_1\circ\gamma_2}a^{v_x^{(1)}v_x^{(2)}}
b^{v_x^{(1)}v_y^{(2)}}c^{v_y^{(1)}v_x^{(2)}}d^{v_y^{(1)}v_y^{(2)}}
r^{v_x^{(2)}}s^{v_y^{(2)}}u)^{-d\chi_g^{(2)}}.$$
The last integral can be decomposed into the product of integrals over
$\{\lambda_2\neq \lambda_1\}$ and $\{\lambda_2=\lambda_1\}$.
The integral over $\{\lambda_2\neq \lambda_1\}$ equals to
$$[\prod_{k<l}(1-(tabcd)^{kv_x^{(1)}}(rs)^{k}u)^{-\mathbb{L}^{-k-l}(\mathbb{L}-1)^2}]^{(\mathbb{L}-2)}.$$
If $\lambda_1=\lambda_2=\lambda$ one can make the affine change
of variables  $A_{\lambda}:(x,y)\mapsto (x,y-\lambda x)$,
and therefore this integral equals to
$${\scriptstyle
\int_{\{v_{\widetilde{y}}^{(2)}>v_x^{(2)}\}}(1-t^{A_{\lambda}(\gamma_1)\circ\gamma_2}(abcd)^{v_x^{(1)}}(rs)^{v_x^{(2)}}u)^{-d\chi_g^{(2)}}=
J_{\sigma^{-1}(A_{\lambda}(\gamma_1))}(t,tabcd,1,1,1,rs\mathbb{L}^{-1},1,u).}$$
The product of the remaining factors is equal to
$${\scriptstyle
\prod_{k<l}(1-(t^{v_x^{(1)}}\cdot (ac)^{v_x^{(1)}}r\mathbb{L}^{-1})^{k}((bd)^{v_x^{(1)}}s\mathbb{L}^{-1})^{l}u)^{-1}\times
\prod_{k>l}(1-((ac)^{v_x^{(1)}}r\mathbb{L}^{-1})^{k}(t^{v_x^{(1)}}\cdot(bd)^{v_x^{(1)}}s\mathbb{L}^{-1}u)^{l})^{-1}\times }$$
$$\times \prod_{k<l}(1-(t^{v_x^{(1)}}\cdot (ac)^{v_x^{(1)}}r\mathbb{L}^{-1}\cdot(bd)^{v_x^{(1)}}s\mathbb{L}^{-1})^{k}\mathbb{L}^{-l}u)^{-(\mathbb{L}-2)}=$$
$$\sum_{k_1,k_2,k_3}\alpha_{k_1,k_2,k_3}(u)t^{k_{1}v_x^{(1)}}(ac)^{k_{2}v_x^{(1)}}(r\mathbb{L}^{-1})^{k_2}(bd)^{k_{3}v_x^{(1)}}(s\mathbb{L}^{-1})^{k_3},$$
therefore
$$\int_{\{v_x^{(1)}=v_y^{(1)}\}}J_{\gamma_1}(t,a,b,c,d,r,s,u)p^{v_x^{(1)}}q^{v_x^{(1)}}d\chi_g^{(1)}=$$
$$\int_{\mathbb{C}^{*}}d\chi_g(\lambda)\int_{\{y^{(1)}=\lambda x^{(1)}+\widetilde{y}^{(1)}\}}
\sum_{k_1,k_2,k_3}\alpha_{k_1,k_2,k_3}(u)(r\mathbb{L}^{-1})^{k_2}(s\mathbb{L}^{-1})^{k_3}\times$$
$$\times J_{\sigma^{-1}(A_{\lambda}(\gamma_1))}(t,tabcd,1,1,1,rs\mathbb{L}^{-1},1,u)(t^{k_1}(ac)^{k_2}(bd)^{k_3}pq)^{v_x^{(1)}}d\chi_g^{(1)}=$$
$$(\mathbb{L}-1)\int_{\{v_{\widetilde{y}}^{(1)}>v_x^{(1)}\}}
\sum_{k_1,k_2,k_3}\alpha_{k_1,k_2,k_3}(u)(r\mathbb{L}^{-1})^{k_2}(s\mathbb{L}^{-1})^{k_3}\times$$
$$\times J_{\sigma^{-1}(\gamma_1)}(t,tabcd,1,1,1,rs\mathbb{L}^{-1},1,u)
(t^{k_1}(ac)^{k_2}(bd)^{k_3}pq)^{v_x^{(1)}}d\chi_g^{(1)}=$$
$$(\mathbb{L}-1)\int_{\mathcal{L}^{(1)}}
\sum\alpha_{k_1,k_2,k_3}(u)(r\mathbb{L}^{-1})^{k_2}(s\mathbb{L}^{-1})^{k_3}\times$$
$$\times J_{\gamma_1}(t,tabcd,1,1,1,rs\mathbb{L}^{-1},1,u)
(t^{k_1}(ac)^{k_2}(bd)^{k_3}pq\mathbb{L}^{-1})^{v_x^{(1)}}d\chi_g^{(1)}=$$
$${\scriptstyle
(\mathbb{L}-1)
\sum\alpha_{k_1,k_2,k_3}(u)(r\mathbb{L}^{-1})^{k_2}(s\mathbb{L}^{-1})^{k_3}
I(t,tabcd,1,1,1,t^{k_1}(ac)^{k_2}(bd)^{k_3}pq,1,rs\mathbb{L}^{-1},1,u).}$$

Combining these integrals, one obtains the following proposition.

\begin{theorem}
$$I(t,a,b,c,d,p,q,r,s,u)=$$
$${\scriptstyle
\sum_{k_1,k_2}\varepsilon_{k_1,k_2}(u)(\mathbb{L}^{-1}r)^{k_1}(\mathbb{L}^{-1}s)^{k_2}
I(t,tabcd,bd,cd,d,(ac)^{k_1}(bdt)^{k_2}pq\mathbb{L}^{-1},c^{k_1}d^{k_2}q,rs\mathbb{L}^{-1},s,u)+}$$
$${\scriptstyle
\sum_{k_1,k_2}\varepsilon_{k_1,k_2}(u)(\mathbb{L}^{-1}s)^{k_1}(\mathbb{L}^{-1}r)^{k_2}
I(t,tabcd,ac,ab,a,(bd)^{k_1}(act)^{k_2}pq\mathbb{L}^{-1},b^{k_1}a^{k_2}p,rs\mathbb{L}^{-1},r,u)+}$$
$${\scriptstyle
(\mathbb{L}-1)
\sum\alpha_{k_1,k_2,k_3}(u)(\mathbb{L}^{-1}r)^{k_2}(\mathbb{L}^{-1}s)^{k_3}
I(t,tabcd,1,1,1,t^{k_1}(ac)^{k_2}(bd)^{k_3}pq,1,rs\mathbb{L}^{-1},1,u).}$$

\end{theorem}

\section*{Acknowledgements}

I would like to thank S. Gusein-Zade for  constant attention,
encouragement and  useful discussions.

Moscow State University,\newline
Independent University of Moscow.\newline
E.mail: gorsky@mccme.ru

\end{document}